 \UseRawInputEncoding
\documentclass[11pt]{article}
\usepackage{amsmath}
\usepackage{mathrsfs}
\usepackage{amsfonts}

\usepackage{amsfonts, amsmath, amssymb}
\usepackage{amssymb,amsfonts,amsmath,color,
latexsym, epsfig,cite, psfrag,eepic,colordvi}
\usepackage{amscd,graphics}
\usepackage{lineno}
\newcommand{\qed}{\hfill $\Box $}
\newcommand{\pf}{\noindent {\bf Proof.} }

\newtheorem{theorem}{Theorem}[section]
\newtheorem{lemma}[theorem]{Lemma}

\newtheorem{conjecture}[theorem]{Conjecture}

\begin{document}

\title{A better  bound on the size of  rainbow matchings}
\author{Hongliang Lu\footnote{Partially supported by the National Natural
Science Foundation of China under grant No.11871391 and
Fundamental Research Funds for the Central Universities}\\
School of Mathematics and Statistics\\
Xi'an Jiaotong University\\
Xi'an, Shaanxi 710049, China\\
\medskip \\
Yan Wang\\
School of Mathematics\\
Georgia Institute of Technology\\
Atlanta, GA 30332, USA\\
\medskip \\
Xingxing Yu\footnote{Partially supported by NSF grant DMS-1600738}\\
School of Mathematics\\
Georgia Institute of Technology\\
Atlanta, GA 30332, USA}

\date{}

\maketitle

\date{}

\maketitle

\begin{abstract}
Aharoni and Howard  conjectured that, for positive integers $n,k,t$
with $n\ge k$ and $n\ge t$,  if $F_1,\ldots, F_t\subseteq {[n]\choose k}$ such
that $|F_i|>{n\choose k}-{n-t+1\choose k}$ for $i\in [t]$ then there exist
$e_i\in F_i$ for $i\in [t]$  such that $e_1,\ldots,e_t$ are
pairwise disjoint. Huang, Loh, and Sudakov proved this conjecture for
$t<n/(3k^2)$. In this paper, we show that this conjecture holds for
$t < n/(2k)$ and $n$ sufficiently large.

%Let $n,k,t$ be three integers such that $k \geq 2$, $n$ is sufficiently large and $n\geq 2kt$. Let $[n]=\{1,2,\ldots, n\}$. In this paper, we show the following: Let $F_1,\ldots, F_t\subseteq {[n]\choose k}$. If $|F_i|>{n\choose k}-{n-t+1\choose k}$ for $1\leq i\leq t$, then there exist
%$e_1\in F_1,\ldots, e_t\in F_t$ such that $e_1,\ldots,e_t$ are pairwise disjoint sets, which gives a linear bound for Rainbow Matching Conjecture.

\end{abstract}

\newpage

\section{Introduction}

For a positive integer $n$, let $[n]$ denote the  set $\{1,\ldots,n\}$.
For a set $S$ with at least $k$ elements, let ${S\choose
  k}=\{e\subseteq S\ : \ |e|=k\}$.
Let $k \ge 2$ be an integer.
A \emph{$k$-uniform hypergraph} or \textit{$k$-graph}  is a pair
$H=(V,E)$, where $V=V(H)$ is a finite set of {\it vertices} and
 $E=E(H)\subseteq {V\choose k}$ is the set of {\it edges}.
We use $e(H)$ to denote the number of edges in $H$.
For any $S \subseteq V(H)$, let $H[S]$ denote the
subgraph of $H$ with $V(H[S])=S$ and $E(H[S])=\{e\in E(H): e \subseteq
S\}$, and let $H-S:= H[V(H)\setminus S]$.

A \emph{matching} in a hypergraph $H$ is a subset of $E(H)$ consisting
of disjoint edges. The maximal size of a matching in a hypergraph $H$
is denoted by $\nu(H)$.  A classical problem in extremal set
theory is to determine $\max e(H)$ with $\nu(H)$ fixed.
Erd\H{o}s
\cite{Erdos65} in 1965 made the following conjecture:
For positive integers $k,n,t$ with $n \ge kt$, every $k$-graph $H$ on $n$ vertices with $\nu(H) <
t$ satisfies $e(H)\leq \max \left\{ {n\choose k}-{n-t+1\choose k}, {kt-1\choose k} \right\}.$
This bound is tight for the complete $k$-graph on $kt-1$ vertices and for the $k$-graph on $n$ vertices in which every edge intersects a fixed set of $t-1$ vertices.
There have been recent activities on this conjecture, see \cite{AHS12,AFH12,FLM,Fr13,Fr17,HLS,LM}.
In particular,  Frankl \cite{Fr13}
proved that if  $n\geq (2t-1)k-(t-1)$ and $\nu(H)<t$ then $e(H)\le
{n\choose k}-{n-t+1\choose k}$, with further improvement  by Frankl
and Kupavskii \cite{FK18}.

There are also attempts to extend the above conjecture of  Erd\H{o}s to a family
of hypergraphs. Let $\mathcal{F} = \{F_1,\ldots, F_t\}$ be a family
of hypergraphs. A set of pairwise disjoint edges, one from
each $F_i$, is called a \emph{rainbow matching} for $\mathcal{F}$. In this case, we also say that ${\cal F}$ or
$\{ F_1,\ldots, F_t\}$ {\it admits} a rainbow matching. Aharoni and
 Howard \cite{AH} made the following conjecture, also see Huang, Loh, and Sudakov \cite{HLS}.
\begin{conjecture}\label{AH-HLS}
Let ${\cal F}=\{F_1,\ldots, F_t\}$ be a  family of subsets in ${[n]\choose k}$. If
\[
e(F_i)> \max\left\{{n\choose k}-{n-t+1\choose k},{kt-1\choose k}\right\}
\]
for all $1\leq i\leq t$, then ${\cal F}$ admits a rainbow matching.
\end{conjecture}

 Huang, Loh, and Sudakov \cite{HLS} proved  that Conjecture \ref{AH-HLS} holds for $n > 3k^2t$.

\begin{theorem}[Huang, Loh, and Sudakov]\label{HLS}
Let $n,k,t$ be three positive integers such that $n > 3k^2t$.
Let ${\cal F}=\{F_1,\ldots, F_t\}$ be a family of subsets of ${[n]\choose k}$. If
\[
e(F_i)> {n\choose k}-{n-t+1\choose k}
\]
for all $1\leq i\leq t$, then ${\cal F}$ has a rainbow matching.
\end{theorem}

Recently, Frankl and Kupavskii \cite{FK20} proved that Conjecture
\ref{AH-HLS} holds when $n\ge 12kt\log(e^2t)$, providing an almost
linear bound. In this paper, we show that Conjecture \ref{AH-HLS} holds when $n > 2kt$ and $n$ is sufficiently large.
 \begin{theorem}\label{main}
%Let $\zeta > 0$ be a real number.
Let $n,k,t$ be three positive integers such that $n > 2kt$ and $n$ is
sufficiently large.
Let ${\cal F}=\{F_1,\ldots, F_t\}$ be a family of subsets of ${[n]\choose k}$.
 If
\[
e(F_i)> {n\choose k}-{n-t+1\choose k}
\]
for all $1\leq i\leq t$, then ${\cal F}$ has  a rainbow matching.
 \end{theorem}

Note that the lower bound on $e(F_i)$ is best possible. Indeed, For $i\in [t]$ let $F_i$ be the $k$-graph on $[n]$ consisting of all edges intersecting
$[t-1]$. Then for $i\in [t]$, $e(F_i)= {n\choose k}-{n-t+1\choose k}$ and $\nu(F_i)=t-1$. Hence, $\{F_1, \ldots, F_t\}$ does not admit any rainbow matching.

This example  naturally corresponds to a special class of
$(k+1)$-graphs ${\cal F}_t(k,n)$. This is defined in Section 2, where
we reduce the problem for finding one such rainbow matching to a problem
about finding ``near'' perfect matchings in a larger  class of
$(k+1)$-graphs, denoted by ${\cal F}^t(k,n)$. This will
allow us to apply various techniques used previously to find large matchings in uniform hypergraphs.

 We show in Section 3 that Theorem~\ref{main} holds when   ${\cal F}^t(k,n)$ is close to ${\cal F}_t(k,n)$, in the sense that most edges of ${\cal F}_t(k,n)$ are also edges of ${\cal F}^t(k,n)$.
To deal with the case   ${\cal F}^t(k,n)$ is not close to ${\cal  F}_t(k,n)$, we follow the approach in \cite{BMS15} and
\cite{ST15}. First, we find a small absorbing matching $M_1$ in ${\cal F}^t(k,n)$ which is done in  Section 4.  (However, the existence of this absorbing matching does not require that ${\cal F}^t(k,n)$ be not close to ${\cal F}_t(k,n)$.) Then we take random samples from  ${\cal F}^t(k,n)-V(M_1)$ so that they satisfy various properties,
in particular they all have fractional perfect matchings, see Section 5.  In Section 6, we use fractional perfect matchings in those random samples to perform a second round of randomization to find a spanning subgraph $H'$ of ${\cal F}^t(k,n)-V(M_1)$. We then
apply a result of Pippenger to find a matching in $H'$ covering all but a small constant fraction of the vertices, and use the matching $M_1$
to find the desired matching in ${\cal F}^t(k,n)$ covering all but fewer than $k$ vertices.

\section{Notation and reduction}

To prove Theorem~\ref{main}, we convert
this rainbow matching problem on $k$-graphs  to a matching problem for a special class of $(k+1)$-graphs.
Let $Q,V$ be two disjoint sets.  A $(k+1)$-graph $H$ with vertex $Q \cup V$
is called \emph{$(1,k)$-partite} with {\it partition classes} $Q,V$ if, for each edge $e\in E(H)$, $|e\cap Q|=1$ and $|e\cap V|=k$.
A $(1,k)$-partite $(k+1)$-graph $H$ with partition classes $Q,V$ is  \emph{balanced} if $|V|=k|Q|$.
We say that $S \subseteq V(H)$ is \textit{balanced} if $|S\cap V|=k|S\cap Q|$.

Let $F_1,\ldots, F_t$ be a family of subsets of  ${[n]\choose k}$ and $X:=\{x_1,\ldots,x_t\}$ be a set of $t$ vertices.
We use $\mathcal{F}^t(k,n)$ to denote the $(1,k)$-partite $(k+1)$-graph  with  partition classes $X,  [n]$ and edge set
\[
E(\mathcal{F}^t(k,n))=\bigcup_{i=1}^t \{\{x_i\}\cup e\ :\ e\in F_i\}.
\]
If $F_1=\cdots=F_t=H_k(t,n)$, where  $H_k(t,n)$ denotes the $k$-graph with vertex set $[n]$
 and edge set  ${[n]\choose k}\setminus {[n]-[t]\choose k}$, then
we denote $\mathcal{F}^t(k,n)$ by $\mathcal{F}_t(k,n)$.

\medskip

\textbf{Observation 1:}  $\{F_1,\ldots,F_t\}$ admits a rainbow matching if, and only if,  $\mathcal{F}^t(k,n)$ has a matching of size $t$.

\medskip

Hence, to prove Theorem~\ref{main}, it suffices to show that $\mathcal{F}^t(k,n)$ has a matching of size $t$. For convenience, we further reduce this problem to a near perfect matching problem.
Write $n-kt=km+r$, where $0\leq r\leq k-1$. Let $F_1, \ldots, F_t\subseteq {[n]\choose k}$, and let $F_i={[n]\choose k}$ for $i=t+1, \ldots, t+m$.
  Let $Q=\{x_1,\ldots,x_{m+t}\}$ and let $\mathcal{H}^t(k,n)$ be the $(1,k)$-partite $(k+1)$-graph with partition classes $Q, [n]$ and edge set
\begin{align*}
E(\mathcal{H}^t(k,n))=\bigcup_{i=1}^{m+t}\{\{x_i\}\cup e\ :\ e\in F_i\}.
\end{align*}
When $F_1=\cdots=F_t=H_k(t,n)$,
we denote $\mathcal{H}^t(k,n)$ by $\mathcal{H}_t(k,n)$. Note that  $\nu(\mathcal{H}_t(k,n))= m+t=(n-r)/k$, i.e., ${\cal H}_t(k,n)$ has a matching covering
all but less than $k$ vertices (and such a matching is said to be near perfect).

\begin{lemma}\label{Rain-PM}
$\mathcal{F}^t(k,n)$ has a  matching of size $t$ if, and only if,  $\mathcal{H}^t(k,n)$ has a matching of size $m+t=\lfloor n/k\rfloor$.
\end{lemma}
\pf First, suppose that  $\mathcal{F}^t(k,n)$ has a matching $M_{1}$ of size $t$. Since $n-kt= km+r\ge km$,  $[n]\setminus V(M_{1})$
contains $m$ pairwise disjoint $k$-sets, say $e_1,\ldots, e_m$. Let $M_{2}=\{e_i\cup \{x_{i+t}\}\ :\ i\in [m]\}$.
Then $M_{1}\cup M_{2}$ is a matching of size $m+t$ in $\mathcal{H}^t(k,n)$.

Now assume that $\mathcal{H}^t(k,n)$ has a matching $M$ of size $m+t$. Note that each edge in $M$ contains exactly one vertex in
$\{x_1, \ldots, x_{m+t}\}$. Thus, the $t$ edges in $M$ containing one of $\{x_1, \ldots, x_t\}$ form a matching in $\mathcal{F}^t(k,n)$ of size $t$. \qed

\medskip

For  the proof of Theorem \ref{main}, we need  additional concepts and notation.
Given two hypergraphs $H_1, H_2$ with $V(H_1)=V(H_2)$, let $c(H_1,H_2)$ be the minimum of $|E(H_1)\backslash E(H')|$
taken over all isomorphic copies $H'$ of $H_2$ with  $V(H') = V(H_2)$.
For a real number $\varepsilon > 0$,
we say that $H_2$ is \textit{$\varepsilon$-close} to $H_1$ if $V(H_1) = V(H_2)$ and $c(H_1,H_2)\leq \varepsilon|E(H_1)|$.
The following is obvious.

\noindent\textbf{Observation 2:} $\mathcal{F}^t(k,n)$ is $\varepsilon$-close to $\mathcal{F}_t(k,n)$ if, and only if,
$\mathcal{H}^t(k,n)$ is $\varepsilon$-close to $\mathcal{H}_t(k,n)$.

%By Theorem \ref{HLS}, we may assume that $n\leq 3k^2 t$.

% using the absorbing method.

%Let $\varepsilon,\varepsilon',\rho,\rho',\gamma,\gamma'$ denote constants such that ***************

As mentioned in Section 1, our proof of Theorem~\ref{main} will be divided into two parts, according to  whether or
not  $\mathcal{F}^t(k,n)$ is $\varepsilon$-close  to $\mathcal{F}_t(n,k)$.
If $\mathcal{F}^t(k,n)$ is close to $\mathcal{F}_t(n,k)$, we will apply greedy argument to construct a matching of size $t$.
If $\mathcal{F}^t(k,n)$ is not close to $\mathcal{F}_t(n,k)$, then, by Observation 2, $\mathcal{H}^t(k,n)$ is not close to $\mathcal{H}_t(n,k)$, and
we will  show that $\mathcal{H}^t(k,n)$ has a spanning subgraph with properties that enable us to find a large matching $M_2$ and to use absorbing matching $M_1$ to enlarge $M_2$ to a near perfect matching.

%Using a R\"{o}dl nibble result, we can find an edge cover of small size.
%Finally, we use an absorbing lemma to construct the desired rainbow matching.

%In Section 2, we deal with the case when  $\mathcal{F}^t(k,n)$ is $\varepsilon$-close to $\mathcal{F}_t(k,n)$ for sufficiently small $\varepsilon$.
%In Section 3, we develop an absorbing lemma which will be used in the subsequent proof. Note that this absorbing lemma does not depend on whether $\mathcal{F}^t(k,n)$ is $\varepsilon$-close to $\mathcal{F}_t(k,n)$.
%In Section 4, we deal with the case  when  $\mathcal{F}^t(k,n)$ is not
%$\varepsilon$-close to $\mathcal{F}_t(k,n)$.
%Let $\varepsilon,\rho,c,\gamma,\xi$ be constants such that $0< \gamma\ll \xi \ll \rho\ll c\ll \varepsilon\ll 1$.

\section{The extremal case: $\mathcal{F}^t(k,n)$ is $\varepsilon$-close to $\mathcal{F}_t(k,n)$}

In this section, we prove Theorem~\ref{main}  for the case when $\mathcal{F}^t(k,n)$
is $\varepsilon$-close to the extremal configuration $\mathcal{F}_t(k,n)$, where, given any real $\zeta$ with $0<\zeta<1$, $\varepsilon$ satisfies
$$2^k \sqrt{\varepsilon} < \min\{((k+1)24^kk^{2k})^{-1},  \zeta^{k-1}(6 k^2 2^k (k-1)!)^{-1}  \}.$$
Note that $\zeta$ will be determined when we consider the non-extremal case where $\mathcal{F}^t(k,n)$ is not
$\varepsilon$-close to $\mathcal{F}_t(k,n)$.

Let $H$ be a $(k+1)$-graph and $v \in V(H)$. We define the neighborhood $N_H(v)$ of $v$ in $H$ to be the set $ \{ S \in {V(H) \choose k} \ :\ S \cup \{v\} \in E(H) \} $.
Let $H$ be a $(k+1)$-graph with the same vertex set as $\mathcal{F}_t(k,n)$.
% and $t = \Theta(n)$.
Given real number $\alpha$ with $0<\alpha < 1$, a vertex $v$ in $H$ is called
\emph{$\alpha$-good} with respect to $\mathcal{F}_t(k,n)$ if
$$\left|N_{\mathcal{F}_t(k,n)}(v)\setminus N_H(v)\right|\le \alpha n^{k} $$
and, otherwise, $v$ is called
\emph{$\alpha$-bad}.
%This notion quantifies the closeness of $H$ to $\mathcal{F}_t(k,n)$ at a vertex.
Clearly, if $H$ is $\varepsilon$-close to $\mathcal{F}_t(k,n)$, then the
number of $\alpha$-bad vertices in $H$ is at most $(k+1)\varepsilon n/\alpha$.

\begin{lemma}\label{good-lem}
Let $\zeta, \alpha $ be real numbers and  $n,k,t$ be positive integers such that $0<\zeta <1$, $n\ge 24k^3$, $t<(1-\zeta)n/k$,
%$n \ge \max(24k^3, 2k/\zeta)$ and $n/(6k^2)\le t<(1-\zeta)n/k$ and
and $\alpha < \min\{((k+1)24^kk^{2k})^{-1},  \zeta^{k-1}(6 k^2 2^k (k-1)!)^{-1}   \}$.
Let $H$ be a $(1,k)$-partite $(k+1)$-graph with $V(H)=V(\mathcal{F}_t(k,n))$.  If every vertex of $H$ is $\alpha$-good  with respect to $\mathcal{F}_t(k,n)$, then $H$ has a matching of size $t$.
\end{lemma}

\pf
Let $X := \{x_1, x_2, ..., x_t\}$, $W := [t]$, and $U := [n]  \setminus [t]$, such that $X,[n]$ are the partition classes of $H$.
Let $M$ be a maximum matching such that $|e\cap X| = |e\cap W|=1$ for all $e\in E(H)$.  Let $X'=X\setminus V(M)$, $W'=W\setminus V(M)$, and $U'=U\setminus V(M)$.

We claim that $|M|\geq n/12k^2$. For,
otherwise, assume $|M| <  n/(12k^2)$.
Consider any vertex $x\in X'$.
Since $x$ is $\alpha$-good, we have
\[
\left| \left(W\times {U\choose k-1}\right)\setminus N_H(x)\right|\leq \alpha n^k.
\]
Note that, since $t<(1-\zeta)n/k$ and $\alpha< \zeta^{k-1}(6k^22^k(k-1)!)^{-1}$,
\[
\left|W'\times {U'\choose k-1}\right|\geq (|W|-|M|){n-kt\choose k-1}>\frac{n}{12k^2}{\zeta n\choose k-1}
>\frac{n}{12k^2} \frac{(\zeta n/2)^{k-1}}{(k-1)!}
> \alpha n^k.
\]
Thus there exists $f\in N_{H}(x)\cap \left(W'\times {U'\choose k-1}\right)$.  Now $M'=M\cup \{\{x\}\cup f\}$ is a  matching of size $|M|+1$ in $H$, and
$|f\cap X|= |f\cap W|=1$ for all $f\in M'$. Hence, $M'$ contradicts the choice of $M$.

Let $S=\{u_1,\ldots,u_{k+1}\}\subseteq V(H)\setminus V(M)$, where $u_1\in X'$, $u_{k+1}\in W'$ and $u_i\in U'$ for $i\in [k]\setminus \{1\}$.
Let $\{e_1,\ldots,e_{k}\}$ be an arbitrary $k$-subset of $M$, and let
$e_i := \{v_{i,1},v_{i,2},\ldots,v_{i,k+1}\}$ with $v_{i,1} \in X$, $v_{i,k+1}\in W$, and $v_{i,j}\in U$ for $i \in [k]$ and  $j \in [k] \setminus \{1\}$.
For $j \in [k+1]$, let $f_j := \{u_{j}, v_{1,j+1},
v_{2,j+2},\ldots,$ $v_{k,j+k}\}$ with addition in the subscripts
modulo $k+1$ (except we write $k+1$ instead of $0$). Note that $f_1, \ldots, f_{k+1}$ are pairwise disjoint.

If $f_j \in E(H)$ for all $j \in [k+1]$ then $M':= (M \cup
\{f_1,\ldots,f_{k+1}\})\setminus \{e_1,\ldots,e_{k}\}$ is a matching in $H$
such that  $|M'| = |M| + 1 > |M|$ and $|f\cap X|=|f\cap W|=1$ for all $f\in M'$, contradicting the choice of $M$.
Hence, $f_j\not\in E(H)$ for some $j \in [k+1]$.

Note that there are $\binom{|M|}{k}k!$  choices of $(e_1,\ldots, e_{k}) \subseteq M^k$ and that for any two different such choices the corresponding $f_j'$s  are distinct.
Hence,
\begin{eqnarray*}
& & |\{e \in E(\mathcal{F}_t(k,n)) \setminus  E(H): |e\cap \{u_i: i\in [k+1]\}|=1\}|\\
&\geq & |M|(|M| - 1) \cdots (|M| - k +1) \\
&> & \left(n/(12k^2) - k \right)^{k} \\
&> & \left(n/(24k^2)\right)^{k} \quad \mbox{ (since $n\ge 24k^3$)}\\
&> & (k+1)  \alpha n^{k}  \quad \mbox{ (since $\alpha< ((k+1)24^{k}k^{2k}))^{-1}$}.
\end{eqnarray*}
This implies that there exists $i \in [k+1]$ such that
$|N_{\mathcal{F}_t(k,n)}(u_i) \setminus N_{H}(u_i)| > \alpha n^{k}$,
contradicting the fact that all $u_i$ are $\alpha$-good. \qed

\medskip

We can now prove Theorem~\ref{main} when $\mathcal{F}^t(k,n)$ is  $\varepsilon$-close to $\mathcal{F}_t(k,n)$.
\begin{lemma}\label{close-lem}
Let $0 < \zeta, \varepsilon < 1$ be real numbers and $k\ge 3$ and $t\ge 0$ be integers, such that
$t< (1-\zeta)(1-k(k+1)\sqrt{\varepsilon}) n/k$, $n\ge 48k^3$,
%$n \ge \max\{ 48k^3, 4k/\zeta\}$
and $2^k \sqrt{\varepsilon} < \min\{((k+1)24^kk^{2k})^{-1},  \zeta^{k-1}(6 k^2 2^k (k-1)!)^{-1}  \}$.
Let $(F_1, \ldots, F_t)$ be a family of subsets of ${[n]\choose k}$ such that
$e(F_i)> {n\choose k}-{n-t+1\choose k}$ for $i\in [t]$,  and  let ${\cal F}^t(k,n)$ denote the corresponding $(1,k)$-partite $(k+1)$-graph.
 Suppose $\mathcal{F}^t(k,n)$ is  $\varepsilon$-close to $\mathcal{F}_t(k,n)$. Then $\mathcal{F}^t(k,n)$ has a matching of size $t$.
\end{lemma}

\pf We may assume $n\le  3k^2t$ as otherwise the assertion follows from Theorem~\ref{HLS}.
Let $B$ denote the set of $\sqrt{\varepsilon}$-bad vertices in $\mathcal{F}^t(k,n)$. Since $\mathcal{F}^t(k,n)$ is
$\varepsilon$-close to $\mathcal{F}_t(k,n)$, $|B|\leq (k+1)\sqrt{\varepsilon}n$.
Let $X, [n]$ be the partition classes of ${\cal F}_t(k,n)$, and let   $X := \{x_1, x_2, ..., x_t\}$, $W := [t]$, and $U := [n] \setminus [t]$.
Note that each edge of ${\cal F}_t(k,n)$ intersects $W$.

Let $b:=\max\{|B\cap X|, |B\cap W|\}$; so $b \le (k+1)\sqrt{\varepsilon}n$.
We choose $X_1,W_1$ such that $B\cap X \subseteq X_1 $, $B\cap W\subseteq W_1$ and $|X_1|=|W_1|=b$.
Let $\mathcal{F}_1=\mathcal{F}^t(k,n)[X_1\cup W_1\cup U]$.
For every $x\in X_1$, we have
\[
|N_{\mathcal{F}_1}(x)|\geq |N_{\mathcal{F}}(x)|-\left({n\choose k}-{n-(t-b)\choose k}\right)
>{n-(t-b)\choose k}-{n-(t-1)\choose k}.
\]
Since $n-(t-b) > n/2 \ge 3k^2(k+1)\sqrt{\varepsilon}n > 3k^2 b$,  it follows from
Theorem \ref{HLS} that the family $\{N_{\mathcal{F}_1}(x)\ |\ x\in X_1\}$  admits a rainbow matching.
%Note that $\mathcal{F}_1$ has a matching of size $b$ if, and only if, $\{N_{\mathcal{F}_1}(x)\ |\ x\in X_1\}$ has a rainbow matching.
Thus, by Observation 1,  $\mathcal{F}_1$ has a matching $M$ of size $b$. Clearly, $M$ covers $B\cap X$.

Let ${\cal F}_2:={\cal F}^t(k,n)[(X\setminus X_1)\cup ([n]\setminus (V(M)\cup B)]$, and let $a:=|B\setminus V(M)|$.
 %Let $\mathcal{F}_2=\mathcal{F}-V(M)-B$ and let $a=|B-V(M)|$.
By the choice of $W_1$ and $X_1$, we have $B\cap (W\setminus W_1)=\emptyset$.   Note that ${\cal F}_2$ may be viewed as    the $(1,k)$-partite $(k+1)$-graph ${\cal F}_2=\mathcal{F}_{t-b}(k, n-kb-a)$, with partiton classes $X\setminus X_1, [n]\setminus (V(M)\cup B)$ from the family
$(F_{i}[(X\setminus X_1)\cup ([n]\setminus (V(M)\cup B)]: i\in X\setminus X_1)$.
Put $n'=n-kb-a$ and $t'=t-b$. We wish to apply Lemma~\ref{good-lem}.

Note that $n'=n-kb-a \ge n - k|B| \ge n-k(k+1)\sqrt{\varepsilon}n \ge n/2 \ge 24k^3$.
%\max(24k^3, 2k/\zeta)$.
 Moreover, since $b \le (k+1) \sqrt{\varepsilon} n \le n/6k^2 \le t/2$, we have $n'/6k^2 \le n/6k^2 \leq t/2 < t-b=t'$.
 Also, $t' \le t  < (1-\zeta)(n-k(k+1)\sqrt{\varepsilon}n)/k  \le (1-\zeta)(n - k|B|)/k  \le (1-\zeta)(n-kb-a)/k= (1-\zeta)n'/k$.

 For every $x\in V(\mathcal{F}_2)$, since $x$ is $\sqrt{\varepsilon}$-good with respect to $\mathcal{F}_t(k,n)$,
 \begin{align*}
|N_{\mathcal{F}_{t'}(k, n')}(x)\setminus N_{\mathcal{F}_2}(x)|&\leq |N_{\mathcal{F}_t(k,n)}(x)\setminus N_{\mathcal{F}}(x)|\\
&\leq \sqrt{\varepsilon} n^k\\
&<2^k\sqrt{\varepsilon}(n-kb-a)^k\quad\mbox{(since $kb+a\leq (k+1)^2\sqrt{\varepsilon}n<n/2$)}\\
&=2^k\sqrt{\varepsilon} (n')^k.
 \end{align*}
Thus every vertex $x$ of $\mathcal{F}_2$ is $2^k\sqrt{\varepsilon}$-good with respect to $\mathcal{F}_{t'}(k,n')$.  By Lemma \ref{good-lem}, $\mathcal{F}_2$ has a matching $M'$ of size $t-b$.
Hence $M\cup M'$ is a matching in $\mathcal{F}$ of size $t$. \qed

%\section{Hypergraphs not-close to $H_t(k,n)$}

\section{Absorbing Lemma}

The purpose of this section is to prove the existence of a small matching $M$ in $\mathcal{H}^t(k,n)$ such that for any small balanced set $S$,
${\cal H}^t(k,n)[V(M)\cup S]$ has a perfect matching. We need to use Chernoff bounds here and in the next section.
Let $Bi(n,p)$ denote a binomial random variable with parameters $n$ and $p$.
The following well-known concentration inequalities, i.e. Chernoff bounds, can be found in Appendix A in \cite{AS08}, or Theorem 2.8, inequalities (2.9) and (2.11) in \cite{JLR}.

 \begin{lemma}[Chernoff inequality for small deviation]\label{chernoff1}
 If $X=\sum_{i=1}^n X_i$, each random variable $X_i$ has Bernoulli distribution with expectation $p_i$, and $\alpha \le 3/2$, then
 $$ \mathbb{P}(|X-\mathbb{E}X| \ge \alpha \mathbb{E}X ) \le 2e^{-\frac{\alpha^2}{3}\mathbb{E}X}.$$
 In particular, when $X \sim Bi(n,p)$ and $\lambda < \frac{3}{2}np$, then
 $$ \mathbb{P}(|X-np| \ge \lambda ) \le e^{-\Omega(\lambda^2/(np))}.$$
 \end{lemma}

% \begin{lemma}[Chernoff inequality for large deviation]\label{chernoff2}
 %If $X=\sum_{i=1}^n X_i$, each random variable $X_i$ has Bernoulli distribution with expectation $p_i$, and $x \ge 7 \mathbb{E}X$, then
 %$$ \mathbb{P}( X \ge x ) \le e^{-x}.$$
 %\end{lemma}

We can now prove an absorbing lemma for $H={\cal H}^t(k,n)$.

\begin{lemma}\label{Absorb-lem}
	Let $k\geq 3$ be an integer, $\zeta > 0$ be a real number an $n\ge n_1(k,\zeta)$ sufficiently large.
%Let $\mathcal{H}=\mathcal{H}^t(k,n)$ be a $(1,k)$-partite $(k+1)$-graph with partition $(X, [n])$ defined above.
Let $H$ be a $(1,k)$-partite $(k+1)$-graph with partition classes $\{x_1, \ldots, x_{\lfloor n/k\rfloor} \}, [n]$ such that
$d_H(x_i)>{n\choose k}-{n-t+1\choose k}$ for $i\in [t]$ and $d_H(x_i)={n\choose k}$ for $i=t+1, \ldots, \lfloor n/k\rfloor$.
Suppose $ n/3k^2 \le t \le (1-\zeta)n/k$.
%Then there exist constant $c(k, \zeta) > 0$ such that, for some $c$ with $0 < c < c(k, \zeta)$,
Then for any $c$ with $0<c<\zeta^{2k} (12k^22^k (k!)^k)^{-2}$, there exists  a matching $M$ in $H$ %$\mathcal{H}$
	such that $|M|\le 2k c n$ and, for any balanced subset $S\subseteq V(H)$ %$S\subseteq V(\mathcal{H})$
       with $|S|\le (k+1)c^{1.5} n/2$,       $H[V(M)\cup S]$ %$\mathcal{H}[V(M)\cup S]$
         has a perfect matching.
\end{lemma}

\pf For balanced $R\in \binom{V(H)}{k+1}$ and balanced $Q\in \binom{V(H)}{k(k+1)}$,
	we say that $Q$ is \emph{$R$-absorbing}
	if  $\nu(H[Q\cup R])=k+1$ and $Q$ is the vertex set of a matching in $H$. Let $\mathcal{L}(R)$ denote the collection of all $R$-absorbing sets in $H$.

\medskip
\textbf{Claim 1.~} For each balanced $(k+1)$-set $R\subseteq V(H)$, %of $\mathcal{H}$,
the number of  $R$-absorbing sets in $H$ is at least $\zeta^k ({n\choose k})^{k+1} (6k^2 2^kk^2!)^{-1}$. %$\mathcal{H}$.

Let $R=\{x,u_1,\ldots,u_k\}$ be fixed with $x\in X$ and $u_i\in [n]$ for $i\in [k]$. Note that the number of edges in $H$ containing
$x$ and intersecting $\{u_1, \ldots, u_k\}$ is at most $k{n\choose k-2}$, and $d_{\mathcal{H}}(x)>{n\choose k}-{n-t+1\choose k}$. So
the number of edges  $\{x,v_1,\ldots, v_k\}$ in $H$ such that
$v_i\in [n]$ for $i\in [k]$ and $\{v_1, \ldots, v_k\}\cap \{u_1,\ldots, u_k\}=\emptyset$ is at least
\begin{align*}
d_{\mathcal{H}}(x)>{n\choose k}-{n-t+1\choose k}-{n\choose k-2}\geq \frac{1}{6k^2}{n\choose k},
\end{align*}
since $3k^2t\geq n\geq kt$.
%i.e., there are at least $\frac{1}{6k^2}{n\choose k} $ edges in $H$ containing $x$.

Fix a choice of an edge $\{x,v_1,\ldots, v_k\}$ in $H$ such that
$v_i\in [n]$ for $i\in [k]$ and $\{v_1, \ldots, v_k\}\cap \{u_1,\ldots, u_k\}=\emptyset$, and
let $W_0 = \{v_1,\ldots, v_k\}$.
For each $j\in [k]$ and each pair $u_j,v_j$, we choose a $k$-set $U_j$ such that $U_j$ is disjoint from
$W_{j-1}\cup R$ and both $U_j\cup \{u_j\}$ and $U_j\cup\{v_j\}$ are edges in $H$, and let $W_{j} :=U_{j}\cup W_{j-1}$.
%For a fixed $j\in [k]$, we call such a choice $U_j$ good.
Then if $W_k$ is defined then  $W_k$ is an absorbing $k(k+1)$-set for $R$.

Note that in each step $j\in [k]$ there are $k+1 +jk$ vertices in $W_{j-1} \cup R$.  Thus, the number of edges in $H$ containing
 $u_j$ (respectively, $v_j$) and at least one other vertex in $W_{j-1} \cup R$ is at most $(k+1 +jk){n\choose k-2} \lfloor n/k \rfloor <(k+1)n{n\choose k-2}$.
Note that by definition of $x_{t+1}, x_{t+2}, ..., x_{\lfloor n/k \rfloor}$, there are at least ${n-2\choose k-1}( \lfloor n/k \rfloor -t) \ge {n-2\choose k-1}\zeta n/k$ sets $U_j$ such that both $U_j\cup \{u_j\}$ and $U_j\cup \{v_j\}$ are edges in $\mathcal{H}$ for large $n$.
Hence, for each $j\in[k]$, there are at least ${n-2\choose k-1}\zeta n/k-(k+1)n{n\choose k-2}\geq \zeta n{n-1\choose k-1}/2k$ such choices for $U_j$ (as $n$ is sufficiently large). Thus, in total we obtain $\frac{1}{6k^2}{n\choose k}(\zeta n{n-1\choose k-1}/2k)^k$ absorbing, ordered $k(k+1)$-sets for $R$, with multiplicity at most $(k^2)!$; so
\[
\mathcal{L}(R)\geq \frac{\frac{1}{6k^2}{n\choose k}(\zeta n{n-1\choose k-1}/2k)^k}{(k^2)!}\geq \frac{\zeta^k {n\choose k}^{k+1}}{6k^2 2^k(k^2)!}.
\]
This completes the proof of Claim 1.

%A $k(k+1)$-subset $S$ of $V(\mathcal{H})$ is called balanced if $|S\cap Y|=k$.
Now, let $c$ be fixed constant with $0<c< \zeta^{2k} (12k^22^k (k!)^k)^{-2}$, and
choose a family ${\cal G}$ of balanced $k(k+1)$-sets by selecting each of the ${ \lfloor n/k \rfloor \choose k}{n\choose k^2}$ balanced sets of size $k(k+1)$ with probability
\[
p := \frac{c n}{{\lfloor n/k \rfloor\choose k}{n\choose k^2}}.
\]
It follows from Lemma~\ref{chernoff1} that, with probability $1-o(1)$, the family ${\cal G}$ satisfies the following properties:
\begin{align}\label{absor-eq1}
|{\cal G}|\leq 2c n
\end{align}
and
\begin{align}\label{absor-eq2}
|\mathcal{L}(R)\cap {\cal G}|\geq p |\mathcal{L}(R)|/2\geq \frac{c \zeta^k n}{12k^2 2^k (k!)^{k}} \geq c^{1.5} n
\end{align}
for all balanced $(k+1)$-sets $R$. Furthermore, we can bound the expected number of intersecting pairs of $k(k+1)$-sets from above by
\begin{align*}
{\lfloor n/k\rfloor \choose k}{n\choose k^2}k(k+1)\left({\lfloor n/k \rfloor-1\choose k-1}{n\choose k^2}+{\lfloor n/k\rfloor \choose k}{n-1\choose k^2-1}\right)p^2\leq c^{1.9} n.
\end{align*}
Thus, using Markov's inequality, we derive that with probability at least $1/2$
\begin{align}\label{absor-eq3}
\mbox{${\cal G}$ contains at most $c^{1.9} n$ intersecting pairs of $k(k+1)$-sets. }
\end{align}

Hence, there exists a family ${\cal G}$ satisfying (\ref{absor-eq1}), (\ref{absor-eq2}) and (\ref{absor-eq3}). Delete one $k(k+1)$-set from each intersecting pair in such a family ${\cal G}$. Further removing all non-absorbing $k(k+1)$-sets, we obtain a subfamily ${\cal G}'$ consisting of pairwise disjoint balanced, absorbing $k(k+1)$-sets, which satisfies
\[
|\mathcal{L}(R)\cap {\cal G}'|\geq \frac{1}{2}c^{1.5} n,
\]
for all balanced $(k+1)$-sets $R$.

Since ${\cal G}'$ consists only of absorbing $k(k+1)$-sets, $H[V ({\cal G}')]$ has a perfect matching $M$,
of size at most $2kcn$ by (\ref{absor-eq2}).  For a balanced set $S\subseteq V(H)$ of size $|S|\leq (k+1)c^{1.5} n/2$, $S$ can be partitioned into at most $c^{1.5} n/2$ balanced $(k+1)$-sets. For each balanced $(k+1)$-set $R$, since $|\mathcal{L}(R)\cap {\cal G}'|\geq \frac{1}{2}c^{1.5} n$,
we can successively choose a distinct absorbing $k(k+1)$-set for $R$ in ${\cal G}'$.
Hence, $\mathcal{H}[V(M)\cup S]$ has a perfect matching. \qed

%\section{$\mathcal{F}^t(k,n)$ is not $\varepsilon$-close to $\mathcal{F}_t(k,n)$}

%In this section, we deal with the case when $\mathcal{F}^t(k,n)$ is not $\varepsilon$-close to $\mathcal{F}_t(k,n)$.

%\subsection{Fractional Perfect Matching}

%A family ${\cal G}$ of subsets    of a set $V$  is
%said to be {\it increasing} if, for any $A\in \mathcal{G}$ and
%$B\subseteq V$,   $A\subseteq B$ implies $B\in \mathcal{G}$.
%Let $H$ be a hypergraph. We  use $\cal{I}(H)$ to
%denote the collection of all independent sets in $H$.

\section{Fractional perfect matchings}

When $\mathcal{F}^t(k,n)$ is not $\varepsilon$-close to $\mathcal{F}_t(k,n)$,
we will use fractional perfect  matchings in random
subgraphs of ${\cal H}^t(k,n)$.

Let $H$ be a hypergraph. A {\it fractional} matching in $H$ is a function $h: E(H) \to [0,1]$ such that
$\sum_{e \ni x} h(e) \le 1$ for all $x \in V(H)$. Let $\nu_f(H):=\max_{h} \sum_{e \in E(H)} h(e)$ which is
the maximum size of fractional matching of $H$. A fractional matching in a $k$-uniform hypergraph with $n$ vertices
is {\it perfect} if its size is $n/k$.

First, we need a concept of dense graphs used in the hypergraph container result of Balogh, Morris, and Samotij  \cite{BMS15} and independently Sexton and Thomassen \cite{ST15}.
Let $H$ be a hypergraph, $\lambda>0$ be a real number, and  ${\cal A}$ be a family of subsets of $V(H)$. We say that $H$ is \textit{$(\mathcal {A}, \lambda)$-dense} if $e(H[A])\ge \lambda e(H)$ for every $A \in \mathcal{A}$.

\begin{lemma}\label{dense}
Let $n, k, t$ be positive integers and $\varepsilon$ be a constant such that $n \le 3k^2t$,
$0<\varepsilon \ll 1$, and $n\geq 40k^2/\varepsilon$. Let $a_0 = \varepsilon/8k, a_1 = \varepsilon/24k^2, a_2=\varepsilon/8k^2$
and $a_3 < \varepsilon/(2^k \cdot k! \cdot 30k)$.
Let $H$ be a $(1,k)$-partite $(k+1)$-graph with vertex partition classes  $X,[n]$ with $|X|=t$.
Suppose $d_H(x)\geq {n\choose k}-{n-t+1 \choose k}- a_3 n^k$ for any $x\in X$. If $H$ is not $\varepsilon$-close to $\mathcal{F}_t(k,n)$, then $H$
is $(\mathcal{A}, a_0)$-dense, where $\mathcal{A}=\{A\subseteq V(H) : |A\cap X|\ge (t/n-a_1) n,\ |A\cap [n]|\ge (1-t/n-a_2) n\}$.
\end{lemma}

\pf We prove this by way of contradiction. Suppose that there exists $A\subseteq V(H)$ such that $|A\cap X|\ge (t/n-a_1) n$,
$|A\cap [n]|\ge  (1-t/n-a_2) n$, and $e(H[A])\le a_0 e(H)$. Without loss of generality, we may choose $A$ such that $|A\cap X|= (t/n-a_1) n$ and
$|A\cap[n]|= (1-t/n-a_2) n$. Let $U\subseteq [n]$ such that $A \cap [n] \subseteq U$ and $|U|=n-t$. Let $A_1=A\cap X$, $A_2=X\setminus A$, $B_1=A \cap [n]$, and $B_2=U\setminus A$.

 Let $H_0$ denote the isomorphic copy of $H$ by naming vertices such that $X = \{x_1, ..., x_t\}$ and $U = [n] \setminus  [t]$.
 We derive a contradiction by showing
that $|E(\mathcal{F}_t(k,n))\setminus E(H_0)|< \varepsilon e(\mathcal{F}_t(k,n))$.
Note that, since $n\le 3k^2t$,
$$e(\mathcal{F}_t(k,n)) = t\left( {n \choose k} - {n-t \choose k}\right)
\geq t\left( {n \choose k} - {n-n/3k^2 \choose k}\right) \ge t {n \choose k}/(3k).$$
Moreover,
$$e(\mathcal{F}_t(k,n)) \ge t {n \choose k}/(3k) = \frac{tn}{3k^2}{ n-1 \choose k-1},$$
and since $n > 2k$,
$$e(\mathcal{F}_t(k,n)) \ge t {n \choose k}/(3k) > \frac{tn^k}{2^k \cdot k! \cdot 3k}.$$

Consider $x \in A$.
Let $E_{H_0}(B_1,x)$ denote the set of edges contained entirely in $B_1 \cup \{x\}$ in $H_0$.
%|\{e | x \in e, e \in E(H_0), V(e) \subseteq B_1 \cup \{x\} \}|
The number of edges in $H_0$ containing $x$ that also exist in $\mathcal{F}_t(k,n)$ is the number of edges in $H_0$ containing $x$ and
 intersecting $[t]$. Hence,
\begin{eqnarray*}
 & &  |\{e : x \in e, e \in E(H_0), e \cap [t] \neq \emptyset\}| \\
&\ge &  d_{H_0}(x) - |\{e : x \in e, e \in E(H_0-[t]), e \cap B_2 \neq \emptyset\}| - |E_{H_0}(B_1,x)| \\
&\ge & \left({n\choose k}-{n-t+1 \choose k}-a_3 n^k\right) - a_2 n {n-t \choose k-1} - |E_{H_0}(B_1,x)|.
\end{eqnarray*}
%By assumption, $d_{H_0}(x)\geq {n\choose k}-{n-t+1 \choose k}-c^{1/4} n^k$.

Therefore, we have
\begin{eqnarray*}
& & |E(\mathcal{F}_t(k,n))\setminus E(H_0)| \\
&= & \sum_{x \in A_1} |\{e : x \in e, e \in E(\mathcal{F}_t(k,n))\setminus E(H_0)\}| + \sum_{x \in A_2} |\{e : x \in e, e \in E(\mathcal{F}_t(k,n))\setminus E(H_0)\}| \\
&\le & \sum_{x \in A_1} \left( {n \choose k} - {n-t \choose k} - |\{e : x \in e, e \in E(H_0), e \cap [t] \neq \emptyset\}|\right) + |A_2| \left( {n \choose k} - {n-t \choose k}\right) \\
&\le & \sum_{x \in A_1} \left[ \left({n \choose k} - {n-t \choose k}\right) - \left({n\choose k} - {n-t+1 \choose k}-a_3 n^k - a_2 n{n-t \choose k-1} - E_{H_0}(B_1,x) \right)\right] \\
& & + a_1 n \cdot e(\mathcal{F}_t(k,n))/t \\
&\le &\sum_{x \in A_1} \left[{n-t+1 \choose k} - {n-t \choose k} + a_3 n^k + a_2 n {n-t \choose k-1} + E_{H_0}(B_1,x)\right] + (3k^2 a_1) \cdot e(\mathcal{F}_t(k,n)) \\
&= &t {n-t \choose k-1} + a_3 t n^k  + a_2 t n {n-t \choose k-1} + \sum_{x \in A_1} E_{H_0}(B_1,x) + (3k^2 a_1) \cdot e(\mathcal{F}_t(k,n)) \\
&\le &(3k^2/n) \cdot e(\mathcal{F}_t(k,n)) + (2^k \cdot k! \cdot 3k a_3) \cdot e(\mathcal{F}_t(k,n)) + (3k^2 a_2) \cdot e(\mathcal{F}_t(k,n)) \\
 & & + e(H_0[A]) + (3k^2 a_1) \cdot e(\mathcal{F}_t(k,n)) \\
&< & a_0 e(H_0) + \left(3k^2/n + 2^k \cdot k! \cdot 3k a_3 + 3k^2 a_2 + 3k^2 a_1\right) \cdot e(\mathcal{F}_t(k,n)) \\
&\le & a_0 t {n \choose k} +\left(3k^2/n + 2^k \cdot k! \cdot 3k a_3 + 3k^2 a_2 + 3k^2 a_1\right) \cdot e(\mathcal{F}_t(k,n)) \\
&\le & \left(3k a_0 + 3k^2/n + 2^k \cdot k! \cdot 3k a_3 + 3k^2 a_2 + 3k^2 a_1\right) \cdot e(\mathcal{F}_t(k,n)) \\
&\leq & \varepsilon \cdot e(\mathcal{F}_t(k,n)),
\end{eqnarray*}
a contradiction since $H$ is not $\varepsilon$-close to $\mathcal{F}_t(k,n)$.  \qed

We  also need a result of Lu, Yu, and Yuan \cite{LYY1}, which is a stability result on matchings in ``stable'' graphs.
For subsets $e=\{u_1, ...,u_k\}, f=\{v_1,...,v_k\} \subseteq [n]$ with $u_i < u_{i+1}$ and $v_i < v_{i+1}$ for $i \in [k-1]$, we write $e \le f$ if $u_i\leq v_i$  for all $i \in [k]$.
A hypergraph $H$ with $V(H) = [n]$ and $E(H) \subseteq {[n] \choose k}$ is said to be \textit{stable} if for $e, f \in {[n] \choose k}$ with $e \le f$, $e \in E(H)$ implies $f \in E(H)$. The following is Lemma 4.2 in \cite{LYY1}.

%In the proof, we need the following stability Lemma.
\begin{lemma}[Lu, Yu and Yuan]\label{stafrankl}
	Let $k$ be a positive integer and let $b$ and $\eta$ be constants, such that $0<b<1/(2k)$ and $0<\eta\le (1+18(k-1)!/b)^{-2}$.
	Let $n,m$ be positive integers such that $n$ is sufficiently large and
	$bn\leq m\leq n/(2k)$.
	Let $H$ be a $k$-graph with vertex set $[n]$.
	Suppose $H$ is stable and  $e(H)>{n\choose k}-{n-m\choose k}-\eta n^{k}$.
	If $H$ is not $\sqrt{\eta}$-close to $H_k(m,n)$, then $\nu(H) > m$.
\end{lemma}

%Let $H$ be a hypergraph.
%A subset $I \subseteq V(H)$ is called an independent set if $H[I]$ has no edges.
%A (fractional) vertex cover $\omega$ of $H$ is a function from $V(H)$ to $[0,1]$ such that
%$\sum_{v \in e} \omega(v) \ge 1$ for all $e \in E(H)$. Let $\nu'_f(H):=\min_{\omega} \sum_{v \in V(H)} \omega(v)$, which is the minimum size of a
%vertex cover.

We now state and prove the main result of this section.

\begin{lemma}\label{Fr-PM}
Let $n, k, t$ be positive integers such that $n\equiv 0\pmod k$ and let $c, \varepsilon$ be constants such that  $0 < c\ll \varepsilon \ll 1$.
Suppose that $n$ is sufficiently large and $n/(3k^2) \le t \le n/(2k)$. Let $H$ be a balanced $(1,k)$-partite $(k+1)$-graph with partition classes $X,[n]$, and let $X=\{x_1,\ldots,x_{n/k}\}$ and $X'=\{x_1,\ldots,x_t\}$.
 Suppose  $d_{H}(x)\geq  {n\choose k}-{n-t+1 \choose k}-\sqrt{c} n^k$  for $x\in X'$, and $d_H(x)={n\choose k}$ for $x\in X\setminus X'$, and assume that
for any independent set $S$ in $H$,  $|S\cap X|\leq (t/n-\varepsilon)n$ or $|S\cap [n]|\leq (1-t/n-\varepsilon)n$.
Then $H$ has a fractional perfect matching.
\end{lemma}

\pf We use linear programming duality between vertex cover and matchings.
Let   $\omega:V(H)\rightarrow [0,1]$ such that $\sum_{v \in e} \omega(v) \ge 1$ for all $e \in E(H)$, and, subject to this,
$\omega(H):=\sum_{v\in V(H)}\omega(v)$ is minimum. (Thus, $\omega$ is a minimum fractional vertex cover of $H$.)
Without loss of generality, we may assume that  $\omega(x_1)\leq \cdots\leq \omega(x_{n/k})$ and $\omega(1)\leq \omega(2)\cdots\leq \omega(n)$.
Let $CL(H)$ be a graph with vertex set $V(H)$ and edge set
\[
E(CL(H))=\left\{e\in {V(H)\choose k+1}\ :\ |e\cap Q|=1\mbox{ and } \sum_{x\in e}\omega(x)\geq 1\right\}.
\]
Note that $H$ is a subgraph of $CL(H)$ and  $\omega$ is also a vertex cover of $CL(H)$. Thus $\omega$ is also a minimum vertex cover of $CL(H)$.

By Linear Programming Duality Theory, we have $\nu_f(H)=w(H)=w(CL(H))=\nu_f(CL(H))$.
% $\nu_f(H)=\nu'_f(H)=\nu'_f(CL(H))=\nu_f(CL(H))$.
Thus it suffices to show that $CL(H)$ has a   fractional perfect matching.
Indeed,  we will prove that $\nu(CL(H))=n/k$, i.e., $CL(H)$ has a perfect matching.

By the definition of $E(CL(H))$, we may assume that
\begin{align}\label{STA-Neigh}
N_{CL(H)}(x_1)\subseteq N_{CL(H)}(x_2)\subseteq\cdots\subseteq N_{CL(H)}(x_{n/k}).
\end{align}
Hence,  $N_{H}(x_i)={[n]\choose k}$ for $i\in [n/k] \setminus [t]$.
It is also easy to see that $N_{H}(x_i)$ is stable for all $i\in [n/k]$.

%We may pick $\eta$ such that $\rho\ll \eta\ll \varepsilon$. We discuss two cases.
%We pick $0 < \eta  \ll \varepsilon \ll 1$ and divide into two cases.
Let $\eta$ be a constant satisfying $c^{1/4}\ll \eta\le \min \{ (1+54k^2(k-1)!)^{-1},  \varepsilon(k(k+1))^{-2}\}$. We distinguish two cases.

\smallskip
\textbf{Case 1.~} $N_{H}(x_1)$ is not $\eta$-close to $H_k(t,n)$.

We observe that $e(N_H(x_1)) = d_H(x_1) \geq  {n\choose k}-{n-t+1 \choose k}-\sqrt{c} n^k
= {n\choose k}-{n-t \choose k} - {n-t \choose k-1} -\sqrt{c} n^k$.
By Lemma \ref{stafrankl} with $m=t$ and $b=1/(3k^2)$, $N_{H}(x_1)$ has a matching $M_1$ of size $t$, and
let $M_1=\{e_1,\ldots,e_t\}$.
By (\ref{STA-Neigh}), $M_1\subseteq N_{CL(H)}(x_i)$ for $i\in [n/k]$. Thus
$M_2=\{e_i\cup \{x_{i}\}\ : \ i\in [t]\}$ is a matching  in $CL(H)$.

Partition $[n]\setminus V(M_2)$ into $n/k-t$ pairwise disjoint $k$-sets, say $f_1,\ldots,f_{n/k-t}$.
Then by  (\ref{STA-Neigh}),
$M_2'=\{f_i\cup \{x_{i+t}\}\ : \ i\in [n/k-t]\}$
is a matching in  $CL(H)\setminus V(M_2)$. Hence $M_2\cup M_2'$ is a perfect matching in $CL(\mathcal{H})$.

\smallskip
\textbf{Case 2.~} $N_{H}(x_1)$ is $\eta$-close to $H_k(t,n)$. (Thus, $N_{CL(H)}(x_1)$ is $\eta$-close to  $H_k(t,n)$.)

Let $B$ denote the set of $\sqrt{\eta}$-bad vertices of $N_{CL(H)}(x_1)$ and let $b=|B|$.
Since $N_{CL(H)}(x_1)$ is $\eta$-close to $H_k(t,n)$, we have $b\leq (k+1)\sqrt{\eta}n$. Consider
$H'=CL(H)- (\{x_{t+1},\ldots,x_{n/k}\}\cup \{n-t+1,\ldots,n\})$.
Note that $kb \le k(k+1)\sqrt{\eta} n < \varepsilon n$; so $b<\varepsilon n/k$.
Since for  any independent set $S$ in $H'$, $|S\cap X|\leq (t/n-\varepsilon)n$ or $|S\cap [n]|\leq (1-t/n-\varepsilon)n$,
we can greedily find pairwise disjoint edges $f_1,\ldots,f_{b}$ in $H'$ such that $x_{t-i+1}\in f_i$ in $H'$. Write $M_{21}=\{f_1,\ldots,f_{b}\}$.

Note that for each vertex $v\in \left([n]\setminus V(M_{21})\right)\setminus B$, we have
\begin{align*}
&|N_{H_k(t-b,n')}(v)\setminus N_{CL(H)-(V(M_{21})\cup  B)}(\{v,x_1\})|\\
\leq &|N_{H_k(t,n)}(v)\setminus N_{CL(H)}(\{v,x_1\})|\\
<&\sqrt{\eta}n^{k-1}\\
<& \eta^{1/3}(n')^{k-1},
\end{align*}
where $n'=|[n]\setminus V(M_{21})\setminus B|$.

Thus, all vertices  of  $N_{CL(H)}(x_1)- (V(M_{21})\cup B)$ in $[n]\setminus V(M_{21})$ are $\eta^{1/3}$-good with respect to $H_k(t-b,n')$. Hence by Lemma \ref{good-lem}, $N_{CL(H)}(x_{1})- (V(M_{21})\cup  B)$ has a matching $M_{22}'$ of size $t-b$. Write $M_{22}'=\{e_1,\ldots, e_{t-b}\}$. By  (\ref{STA-Neigh}),
$M_{22}=\{e_i\cup \{x_{i}\}\ :\ i \in [t-b]\}$ is a matching in $H'$.
Thus, $M_{22}\cup M_{21}$ is a matching of size $t$ in $H'$.

Partition $[n]\setminus V(M_{21}\cup M_{22})$ into $n/k-t$ disjoint $k$-sets, say $g_1,\ldots, g_{n/k-t}$. Let
$M_{23}=\{g_i\cup \{x_{i+t}\}\ :\ i\in [n/k]\setminus [t]\}$. Then $M_{21}\cup M_{22}\cup M_{23}$ is a perfect matching in $CL(H)$.
This competes the proof. \qed

\section{Random rounding}

In this section, we will complete the proof of Theorem~\ref{main}. For convenience, in this section we will not round certain numbers to integers this does not
affect calculations.

  First, we need another result of
Lu, Yu, and Yuan \cite{LYY2} on  the independence number of
a subgraph of a $k$-graph induced by a random subset of vertices, which is a generalization of Lemma 4.3 in \cite{LYY2} where it was shown for
$(1,3)$-partite graphs. The same proof for Lemma 4.3 in \cite{LYY2} works here as well  by using Lemma~\ref{dense} in the place of Lemma 4.1 in \cite{LYY2}.

% Note that
%whose edges are well distributed.
\begin{lemma}[Lu, Yu, and Yuan]\label{indep}
        Let $l, \varepsilon', \alpha_1,\alpha_2$ be positive reals, let $\alpha>0$ with  $\alpha \ll \min\{\alpha_1,\alpha_2\}$,
       let $k,n$ be positive integers, and let
       $H$ be a $(1,k)$-partite $(k+1)$-graph  with partition classes $Q,P$ such that $k|Q|=|P|=n$, $e(H)\ge ln^{k+1}$,  and $e(H[F])\ge
        \varepsilon' e(H)$ for all $F\subseteq
        V(H)$ with $|F\cap P|\ge \alpha_1 n$ and $|F\cap Q|\ge \alpha_2 n$.
         Let $R\subseteq V(H)$ be obtained  by taking each vertex of
           $H$ uniformly at random with probability $n^{-0.9}$.
        Then, with probability at least $1-n^{O(1)}e^{-\Omega (n^{0.1})}$, every independent set $J$ in $H[R]$
             satisfies $|J\cap P|\le (\alpha_1 +\alpha+o(1))n^{0.1}$ or $|J\cap Q|\le (\alpha_2 +\alpha+o(1))n^{0.1}$.
\end{lemma}

Next, we also need the Janson's inequality to provide an exponential upper bound for the lower tail of a sum of dependent zero-one random variable. (See Theorem 8.7.2 in \cite{AS08})

\begin{lemma}[Janson] \label{janson}
Let $\Gamma$ be a finite set and $p_i \in [0,1]$ be a real for $i \in \Gamma$.
Let $\Gamma_{p}$ be a random subset of $\Gamma$ such that the elements are chosen independently with $\mathbb{P}[i \in \Gamma_p] = p_i$ for $i \in \Gamma$.
Let $S$ be a family of subsets of $\Gamma$.
For every $A \in S$, let $I_A = 1$ if $A \subseteq \Gamma_p$ and $0$ otherwise.
Define $X = \sum_{A \in S} I_A$,
$\lambda = \mathbb{E}[X]$,
$\Delta = \frac{1}{2}\sum_{A \neq B} \sum_{A \cap B \neq \emptyset} \mathbb{E}[I_A I_B]$
and $\bar{\Delta} = \lambda + 2\Delta$.
Then, for $0 \le t \le \lambda$, we have
$$\mathbb{P}[X \le \lambda - t] \le \exp(-\frac{t^2}{2\bar{\Delta}}). $$
\end{lemma}

Now,  we use Chernoff bound and Janson's inequality to prove a result on several properties of  certain random   subgraphs.

 \begin{lemma}\label{lem1-5}
        Let $n, k$ be integers such that $n\ge k\geq 3$,
       let $H$ be a $(1,k)$-partite $(k+1)$-graph with partition classes $A,B$
       and  $k|A| = |B| = n$,
       let $A_1,A_2$ be a partition of $A$ with $|A_1|\ge n/(3k^2)$ and $|A_2|\ge n/(3k^2)$, and
       let $A_3\subseteq A$ and $A_4\subseteq B$ with $|A_i|=n^{0.99}$ for $i=3,4$.
	Take $n^{1.1}$ independent copies of $R$ and denote them by $R^i$, $1\le i\le n^{1.1}$, where $R$ is chosen from $V(H)$ by taking each vertex uniformly at random with probability $n^{-0.9}$ and then deleting $O(n^{0.06})$ vertices  uniformly at random so that $|R|\in (k+1) \mathbb{Z}$ and $k|R\cap A|=|R\cap B| $.
        For each $S\subseteq V(H)$, let $Y_S:=|\{i: \ S\subseteq R^i\}|$.
            Then, with probability at least $1-o(1)$,  all of the following statements hold:
        \begin{itemize}
            \item [$(i)$] $Y_{\{v\}}=(1 \pm n^{-0.01}) n^{0.2}$  for all  $v\in V(H)$.  %\sim n^{0.2},
            \item [$(ii)$] $Y_{\{u,v\}}\le 2$ for all $\{u, v\} \subseteq V(H)$.
            \item [$(iii)$] $Y_e\le 1$ for all  $e \in E(H)$.
            \item [$(iv)$] For all $i= 1, \dots ,n ^{1.1}$, we have
              $|R_i\cap A| =(1/k\pm o(n^{-0.04}))n^{0.1}$ and  $|R_i\cap B| =(1\pm o(n^{-0.04}))n^{0.1}$,
            \item [$(v)$] Suppose $n/k^3 \le m\le n/k$ and  $\rho$ is a constant with $0<\rho <1$ such that
             $d_{H}(v)\ge {n\choose
                   k}-{n-m\choose k}-\rho n^{k}$ for all $v\in A$.  Then
                   %there exists $\rho' > 0$ such that
                for $1 \le i \le n^{1.1}$ and $v\in R_i \cap A$, we have
               $$d_{R_i}(v)>  {|R_i \cap B|\choose k}-{|R_i \cap B|-mn^{-0.9} \choose k}-3\rho |R_i \cap B|^{k},$$
           \item[$(vi)$] $|R_i\cap A_j|= |A_j|n^{-0.9}\pm n^{0.06}$ for $1 \le i \le n^{1.1}$ and $j\in \{1,2,3,4\}$.
            %\item[$(vii)$] for any idenpendent set $S$ of $R_i$, $|S\cap R_i\cap Q|< t(1-\varepsilon/2)$ and $|S\cap [n]\cap R_i|< (n-t)(1-\varepsilon)$
%with probability at most $n^{O(1)}e^{-\Omega (n^{0.1})}$.
        \end{itemize}
    \end{lemma}

\pf
For $1 \le i \le n^{1.1}$ and $j\in \{1,2,3,4\}$, $\mathbb{E}[|R_i \cap A|] = n^{0.1}/k$, $\mathbb{E}[|R_i \cap B|] = n^{0.1}$ and $\mathbb{E}[|R_i \cap A_j|] = n^{-0.9}|A_j|$. Recall the assumptions $|A_1|\ge n/(3k^2)$,  $|A_2|\ge n/(3k^2)$, and $|A_3|=|A_4|=n^{0.99}$.
By Lemma~\ref{chernoff1}, we have
 \begin{itemize}
 \item []  $\mathbb{P}\left(\left||R_i \cap A| - n^{0.1}/k\right| \ge n^{0.06} \right) \le e^{-\Omega(n^{0.02})}$,
 \item [] $ \mathbb{P}\left(\left||R_i \cap B| - n^{0.1}\right| \ge n^{0.06} \right) \le e^{-\Omega(n^{0.02})}$, and
  \item [] $\mathbb{P}\left(\left||R_i \cap A_j| - |A_j|n^{-0.9}\right| \ge n^{0.06} \right) \le e^{-\Omega(n^{0.02})}$.
\end{itemize}
 Hence, with probability at least $1-O(n^{1.1})e^{-\Omega(n^{0.02})}$, $(iv)$ and $(vi)$ hold.

 For every $v\in V(H)$, $\mathbb{E}[Y_{\{v\}}]=n^{1.1} \cdot n^{-0.9}= n^{0.2}$. By Lemma~\ref{chernoff1},
 $$ \mathbb{P}\left(\left||Y_{\{v\}}| - n^{0.2} \right| \ge n^{0.19} \right) \le e^{-\Omega(n^{0.18})}$$
 Hence, with probability at least $1-O(n)e^{-\Omega(n^{0.18})}$, $(i)$ holds.

 Let $Z_{p,q} = \left|S \in  {V(H) \choose p} : Y_S \ge q \right|$. Then
 $$\mathbb{E}\left[Z_{p,q}\right] \le  {n \choose p} {n^{1.1} \choose q} (n^{-0.9})^{pq} \le n^{p + 1.1q - 0.9pq}. $$
 So $\mathbb{E}[Z_{2,3}] \le n^{-0.1}$ and $\mathbb{E}Z_{k,2} \le n^{2.2 - 0.8k} \le n^{-0.2}$ for $k \ge 3$.
 % by Lemma~\ref{chernoff2},
Hence by Markov's inequality, $(ii)$ and $(iii)$ hold  with probability at least $1-o(1)$.

 Finally we show $(v)$.
 For all $v\in A$, since
             $d_{H}(v)\ge {n\choose
                   k}-{n-m\choose k}-\rho n^{k}$,  we see that, for $1 \le i \le n^{1.1}$ and $v\in R_i \cap A$,
               $$\mathbb{E}\left[d_{R_i}(v)\right]> {n\choose k} n^{-0.9k}-{n-m\choose k}n^{-0.9k}-\rho n^{k}n^{-0.9k}
               > {n^{0.1} \choose k} -{n^{0.1} - mn^{-0.9} \choose k}-\rho n^{0.1k}.$$
By $(iv)$, with probability at least $1-O(n^{1.1})e^{-\Omega(n^{0.02})}$, for all $i= 1, \dots ,n ^{1.1}$, we have
              $|R_i\cap B| =(1+o(n^{-0.04}))n^{0.1}$.
Thus for all $v\in A\cap R_i$,
%there exists some $\rho'' >0$ such that
               $$\mathbb{E}\left[d_{R_i}(v)\right] > {|R_i \cap B|\choose k}-{|R_i \cap B|-mn^{-0.9} \choose k}-2\rho |R_i \cap B|^{k}.$$
We wish to apply Lemma~\ref{janson} with $\Gamma = B$, $\Gamma_p = R_i$ and
$S$ be a family of all $k$-set of $B$.
We define
$$\Delta = \frac{1}{2} \sum_{b_1, b_2 \subseteq B, b_1 \ne b_2, b_1 \cap b_2 \ne \emptyset} \mathbb{E}[I_{b_1}I_{b_2}] \le \frac{1}{2} |R_i \cap B|^{2k-1} $$
By Lemma~\ref{janson},
\begin{align*}
&\mathbb{P}\left( d_{R_i}(v) \leq {|R_i \cap B|\choose k}-{|R_i \cap B|-mn^{-0.9} \choose k}-3\rho |R_i \cap B|^{k} \right) \\
\leq & \mathbb{P} \left( d_{R_i}(v) \leq \mathbb{E}[d_{R_i}(v)] - \rho |R_i \cap B|^{k} \right) \\
\leq & \exp(- \frac{\rho^2 |R_i \cap B|^{2k}}{2{|R_i \cap B|\choose k} + 2|R_i \cap B|^{2k-1}}) \\
\leq & \exp(-\Omega(n^{0.1})).
\end{align*}
Therefore, with probability at least $1-O(n^{1.1})e^{-\Omega(n^{0.1})}$, $(v)$ holds.

By applying  union bound, $(i)$ -- $(v)$ all hold with probability $1-o(1)$.
\qed

\medskip

%\subsection{The Proof of Theorem \ref{main}}

Now we use random subgraphs and fractional matchings to perform a second round of randomization to find a sparse subgraph
in a hypergraph that is not $\varepsilon$-close to $\mathcal{H}_t(k,n)$.

\begin{lemma}\label{Span-subgraph}
Let $k\ge 3$ be an integer, $0 < \rho \ll \varepsilon \ll 1$ be reals, and $n\in k\mathbb{Z}$ be sufficiently large.
Suppose $n/(3k^2) \le t \le n/(2k)$.
Let $H$ is a $(1,k)$-partite $(k+1)$-graph with partition classes $A,B$ such that $k|A|=|B|=n$.
Let  $A_1$ and $A_2$ be a partition of $A$ such that $|A_1|=t$ and $|A_2|=n/k-t$. Suppose that $d_{H}(x)>{n\choose k}-{n-t+1\choose k}-\rho n^k$ for all $x\in A_1$ and $d_{H}(x)={n\choose k}$ for all $x\in A_2$. If $H$ is not $\varepsilon$-close to $\mathcal{H}_t(k,n)$, then there exists a spanning subgraph $H'$ of $H$
such that the following conditions hold:
\begin{itemize}
		\item[$(1)$] For all $x\in V(H')$, with at most $n^{0.99}$ exceptions,
                               $d_{H'}(x)=(1\pm n^{-0.01})n^{0.2}$;
		\item[$(2)$] For all $x\in V(H')$, $d_{H'}(x)< 2 n^{0.2}$;
		\item[$(3)$] For any two distinct $x,y\in V(H')$, $d_{H'}(\{x,y\})< n^{0.19}$.
	\end{itemize}
\end{lemma}

\pf
Let $A_3\subseteq A$ and $A_4\subseteq B$ with $|A_i|=n^{0.99}$ for $i=3,4$.
Let $R_1,\ldots,R_{n^{1.1}}$ be defined as in Lemma \ref{lem1-5}. By Lemma \ref{lem1-5} $(iv)$, we have, for all $i= 1, \dots ,n ^{1.1}$,
              $$|R_i\cap A| =(1/k+o(n^{-0.04}))n^{0.1} \mbox{ and } |R_i\cap B| =(1+o(n^{-0.04}))n^{0.1}.$$
By Lemma \ref{lem1-5} $(vi)$, we have
			  $$|R_i\cap A_1| =(t/n+o(n^{-0.04}))n^{0.1} \mbox{ and } |R_i\cap A_2| =(1/k-t/n+o(n^{-0.04}))n^{0.1}.$$
By Lemma \ref{lem1-5} $(v)$, we have  for $1 \le i \le n^{1.1}$ and $x\in A \cap R^i$,
               $$d_{R_i}(x)>  {|R_i \cap B|\choose k}-{|R_i \cap B|-(t-1)n^{-0.9} \choose k}-3\rho |R_i \cap B|^{k};$$

By $(iv)$ and $(vi)$ of Lemma \ref{lem1-5},
we may choose $I_i\subseteq R_i\cap (A_3\cup A_4)$ such that $R_i\setminus I_i$ is balanced and  $|R_i'|=(1 - o(1))|R_i|$, where $R_i'=R_i\setminus I_i$ for $i=1, \ldots, n^{1.1}.$
Let $H_1=H[A_1 \cup B]$.

Since $H$ is not $\varepsilon$-close to $\mathcal{H}_t(k,n)$, $H_1$ is not $\varepsilon$-close to $\mathcal{F}_t(k,n)$ by Observation 2 in Section 2.
Let $a_0 = \varepsilon/(8k), a_1 = \varepsilon/(24k^2), a_2=\varepsilon/(8k^2)$, and $a_3 < \varepsilon(2^k \cdot k! \cdot 30k)^{-1}$.
By  applying Lemma \ref{dense} to $H_1, a_0, a_1, a_2, a_3$,  we see that $H_1$ is $(\mathcal{F},a_0)$-dense, where
$$\mathcal{F}=\{U\subseteq V(H) : |U\cap A_1|\ge (t/n-a_1 ) n,\ |U\cap B|\ge (1-t/n-a_2 ) n\}.$$

%{\color{red}Delete ``Note that $e(H_1) \ge \beta e(H)$ where $\beta = (1 - (1-1/3k^2)^k)/3k \ge 1/9k^3$.
%Hence, $H$ is $(\mathcal{F},\varepsilon/54k^4)$-dense."}

Now we apply Lemma \ref{indep} to $H_1$ with $l = (3k^3k!)^{-1}$, $\alpha_1 = t/n - a_1$, $\alpha_2 = 1-t/n-a_2$, and $\varepsilon' =a_0$.
Therefore, with probability at least $1-n^{O(1)}e^{-\Omega (n^{0.1})}$, for any independent set $S$ of $R_i'$,  $|S\cap R_i'\cap A_1|\leq (t/n - a_1+o(1))n^{0.1}$ or $|S\cap R_i'\cap B|\leq (1-t/n-a_2+o(1))n^{0.1}$.
By definition, for $x \in R_i' \cap A_2$, $d_{R_i'}(x) = {|R_i'| \choose k}$.

By applying Lemma \ref{Fr-PM} to each $H[R_i']$, we see that each $H[R_i']$ contains a fractional perfect matching $\omega_i$.
Let $H^*=\cup_{i=1}^{n^{1.1}}R_i'$. We select a generalized binomial subgraph $H'$ of $H^*$ by letting $V(H')=V(H)$ and
independently choosing edge $e$ from $E(H^*)$,
with probability $\omega_{i_e}(e)$ if $e\subseteq R_{i_e}'$. (By  Lemma \ref{lem1-5} $(iii)$, for each $e\in E(H^*)$,  $i_e$ is uniquely defined.)

Note that since $w_i$ is a fractional perfect matching of $H[R_i']$ for $1 \le i \le n^{1.1}$, $\sum_{e \ni v} w_i(e) \le 1$ for $v \in R_i'.$
By  Lemma \ref{lem1-5} $(i)$ and by Lemma~\ref{chernoff1},   $d_{H'}(v)=(1 \pm n^{-0.01}) n^{0.2}$ for any vertex $v\in V(H)-(\cup_{i=1}^{n^{1.1}} I_i) \subseteq V(H)-(A_3 \cup A_4)$
and $d_{H'}(v) \le (1 \pm n^{-0.01}) n^{0.2} < 2n^{0.2}$ for vertex $v \in \cup_{i=1}^{n^{1.1}} I_i$.
By Lemma \ref{lem1-5} (ii) $d_{H'}(\{x,y\})\leq 2 < n^{0.19}$ for any $\{x,y\}\in {V(H) \choose 2}$.
%We choose $D=n^{0.2}$, $\tau=\tau'=n^{-0.01}$. Then
Therefore, $H'$ is the desired hypergraph.
\qed

To prove Theorem~\ref{main}, we also need the following result which was  attributed to Pippenger \cite{PS}
(see Theorem 4.7.1 in \cite{AS08}). An {\it edge cover} in a hypergraph $H$
is a set of edges whose union is $V(H)$.
%Then we use a theorem due to Pippenger \cite{PS}, following Frankl and R\"odl \cite{FR} (also see Theorem 8.4 in \cite{Fu}, Theorem 4.7.1 in \cite{AS}).
\begin{theorem}[Pippenger]\label{nibble}
	For every integer $k\ge 2$ and real $r\ge 1$ and $a>0$, there are $\gamma=\gamma(k,r,a)>0$ and $d_0=d_0(k,r,a)$ such that for every $n$ and $D\ge d_0$ the following holds: Every $k$-uniform hypergraph $H=(V,E)$ on a set $V$ of $n$ vertices in which all vertices have positive degrees and which satisfies the following conditions:
	\begin{itemize}
	 	\setlength{\itemsep}{0pt}
		\setlength{\parsep}{0pt}
		\setlength{\parskip}{0pt}
		\item[$(1)$] For all vertices $x\in V$ but at most $\gamma n$ of them, $d_H(x)=(1\pm \gamma)D$;
		\item[$(2)$] For all $x\in V$, $d_H(x)<r D$;
		\item[$(3)$] For any two distinct $x,y\in V$, $d_H(\{x,y\})<\gamma D$;
	\end{itemize}
	contains an edge cover of at most $(1+a)(n/k)$ edges.
\end{theorem}

\noindent\textbf{Proof of Theorem \ref{main}.}
By Theorem \ref{HLS}, we may assume that $2kt < n\leq 3k^2 t$.
Let $0 < \varepsilon \ll 1$ be sufficiently small and $n$ be sufficiently large.
By Observation 1, it suffices to show $\mathcal{F}^t(k,n)$ has a matching of size $t$.
Applying Lemma \ref{close-lem} to $\mathcal{F}^t(k,n)$ with $\zeta = 1/3$,
we may assume that $\mathcal{F}^t(k,n)$ is not $\varepsilon$-close to $\mathcal{F}_t(k,n)$. That is,
$\mathcal{H}^t(k,n)$ is not $\varepsilon$-close to $\mathcal{H}_t(k,n)$ by Observation 2.

Now we apply Lemma \ref{Absorb-lem} to $\mathcal{H}^t(k,n)$ with $\zeta = 1/2$.
Thus there exists some constant $0 < c \ll \varepsilon$ such that $n-kcn \ge 2kt$ and $\mathcal{H}^t(k,n)$ contains an absorbing matching $M_1$ with $m_1:=|M_1|\leq c n$ and for any balanced subset $S$ of vertices with $|S|\leq (k+1)c^{1.5} n$, $\mathcal{H}^t(k,n)[V(M_1)\cup S]$ has a perfect matching. Let $H:=\mathcal{H}^t(k,n)-V(M_1)$ and $n' := n-km_1$.

Next, we see that $H$ is not $(\varepsilon/2)$-close to $\mathcal{H}_{t}(k,n-km_1)$. For, suppose otherwise.
Then
\begin{align*}
&|E(\mathcal{H}_t(k,n)) \setminus  E(\mathcal{H}^t(k,n))| \\
&\le |E(\mathcal{H}_{t}(k,n-km_1)) - E(H)| + |e \in E(\mathcal{H}_t(k,n)) : e \cap V(M_1) \neq \emptyset| \\
&\leq (\varepsilon/2) |E(\mathcal{H}_{t}(k,n-km_1))| + (k+1)cn \cdot n^k \\
&\le \varepsilon |E(\mathcal{H}_t(k,n))|.
\end{align*}
This is a contradiction as   $\mathcal{H}^t(k,n)$ is not $\varepsilon$-close to $\mathcal{H}_t(k,n)$.

Since $n' \ge n - kcn \ge 2kt$, by Lemma \ref{Span-subgraph} $H$ has a spanning subgraph $H'$ such that
\begin{itemize}
		\item[(1)] For all vertices $x\in V(H')$ but at most $n'^{0.99}$ of them,
                $d_{H'}(x)=(1\pm n'^{-0.01})n'^{0.2}$;
		\item[(2)] For all $x\in V(H')$, $d_{H'}(x)< 2 n'^{0.2}$;
		\item[(3)] For any two distinct $x,y\in V(H')$, $d_{H'}(\{x,y\})< n'^{0.19}$.
	\end{itemize}

Hence by applying Lemma \ref{nibble} to $H'$ with $0 < a \ll c^{1.5}$,
$H'$ contains an edge cover of at most $(1+a)((n'/k+n')/(k+1))$ edges.
Thus, at most $a(n'/k+n')$ vertices are each covered by more than one edge
in the cover. Hence,
after removing at most $a(n'/k+n')$ edges from the edge cover, we obtain a matching $M_2$ covering all but at most $(k+1)a(n'/k+n') \le 3kan' \le 3kan$ vertices.

Now we may choose a balanced subset $S$ of $V(H)\setminus V(M_2)$ such that $|V(H)\setminus (V(M_2)\cup S)|\leq k$.
Since $|S| \le 3kan \le (k+1)c^{1.5}n$, $\mathcal{H}^t(k,n)[V(M_1)\cup S]$ has a perfect matching, say $M_3$. Thus,
$M_2\cup M_3$ is matching of $H^t(k,n)$ covering all but at most $k$ vertices, and, hence, has size $\lfloor n/k\rfloor$.
Therefore, by Lemma \ref{Rain-PM}, $\mathcal{F}^t(k,n)$ has a matching of size $t$. \qed


\begin{thebibliography}{99}
%\bibitem{AH09} R. Aharoni and E. Berger, Rainbow matchings in $r$-partite $r$-graphs, \emph{Electron. J. Combin.}, \textbf{16} (2009), \#R119.

%\bibitem{AH17}  R. Aharoni and D. Howard, A rainbow $r$-partite version of the Erd\"os-Ko-Rado Theorem, \emph{Combinatorics, Probability and Computing}, \textbf{26} (2017), 321--337.

\bibitem{AH}  R. Aharoni and D. Howard, Size conditions for the
existence of rainbow matching, Preprint.


\bibitem{AHS12} N. Alon, H. Huang, and B. Sudakov, Nonnegative $k$-sums, fractional
covers, and probability of small deviations, {\it J. Combin. Theory
  Ser. B},  {\bf 102} (2012), 784--796.

\bibitem{AFH12} N. Alon, P. Frankl, H. Huang, V. R\"{o}dl, A. Ruci\'{n}ski, and
  B. Sudakov, Large matchings in uniform hypergraphs and the conjectures of
Erd\H{o}s and Samuels, {\it  J. Combin. Theory Ser. A},  {\bf 119} (2012), 1200--1215.

\bibitem{AS08} N. Alon, J. Spencer, The Probabilistic Method, John Wiley, Inc., New York, 2008.

\bibitem{BMS15}
J. Balogh, R. Morris and W. Samotij, Independent sets in hypergraphs,
{\it  J. American Math. Soc.}, {\bf 28} No. 3 (2015), 669--709.



%\bibitem{BDE76} B. Bollob\'as, D. Daykin and P. Erd\H{o}s, Sets of independent edges of a hypergraphs, \emph{Quart. J. Math. Oxford Ser.}, \textbf{27} (1976), 25--32.

\bibitem{Erdos65}P. Erd\H{o}s,   A problem on independent $r$-tuples, \emph{Ann. Univ. Sci. Budapest. E\"otv\"os Sect. Math.}, \textbf{8} (1965),
93--95.

\bibitem{FLM}
P. Frankl, T. {\L}uczak, K. Mieczkowska, On matchings in hypergraphs, \emph{Electron. J. Combin.}, \textbf{19} (2012), \#R42.

\bibitem{Fr13} P. Frankl, Improved bounds for Erd\H{o}s'
 matching conjecture, {\it J. Combin. Theory Ser. A}, {\bf 120} (2013), 1068--1072.

\bibitem{Fr17}
P. Frankl, On the maximum number of edges in a hypergraph with given matching number, \emph{Discrete Appl. Math.}, \textbf{216} (2017), 562--581.

\bibitem {FK18} P. Frankl and A. Kupavskii, The  Erd\H{o}s matching
  conjecture and concentration inequalities, arXiv:1806.08855v2
  [math.CO].

%{\color{blue}
\bibitem {FK20} P. Frankl and A. Kupavskii, Simple juntas for shifted families, \emph{Discrete Analysis}, (2020), 14507.
%}

\bibitem{HLS} H. Huang, P. Loh, and B. Sudakov, The size of a hypergraph and its matching number, \emph{Combinatorics, Probability and Computing,} \textbf{21}  (2012), 442--450.


\bibitem{JLR} S. Janson, T. {\L}uczak, A. Ruci\'{n}ski, Random Graphs, John Wiley and Sons, New York, 2000.

\bibitem{LYY1} H. Lu, X. Yu and X. Yuan, Nearly perfect matchings in uniform hypergraphs, arXiv:1911.07431 [math.CO].

\bibitem{LYY2}
H. Lu, X. Yu and X. Yuan, Rainbow matchings for 3-uniform hypergraphs, arXiv:2004.12558 [math.CO].

\bibitem{LM}
    T.  {\L}uczak and K. Mieczkowska. On Erd\"os¡¯ extremal problem on matchings in hypergraphs, \emph{J. Combin. Theory Ser. A},
\textbf{124} (2014), 178--194.

%\bibitem{MN}
%M. Matsumoto and N. Tokushige, The exact bound in the Erd\H{o}s-Ko-Rado theorem for cross-intersecting families, \emph{J. Combin.
%Theory Ser. A}, \textbf{52} (1989), 90--97.

\bibitem{PS}
N. Pippenger and J. Spencer, Asymptotic behaviour of the chromatic index for hypergraphs, \emph{J. Combin. Theory, Ser. A}, \textbf{51} (1989), 24--42.


%\bibitem{PY} L. Pyber, A new generalization of the Erd\H{o}s-Ko-Rado theorem, \emph{J. Combin. Theory
%Ser. A}, \textbf{43} (1986), 85--90.


\bibitem{ST15} D. Saxton and A. Thomason, Hypergraph containers, {\it
    Invent. Math.}, {\bf 201} (3) (2015) 925--992.
%
%
%
\end{thebibliography}
\end{document}